\documentstyle[amscd,amssymb,verbatim,epsf,mathrsfs,graphicx,hyperref]{amsart}

\begin{document}
\theoremstyle{plain}
\newtheorem{Thm}{Theorem}
\newtheorem{Cor}{Corollary}
\newtheorem{Con}{Conjecture}
\newtheorem{Main}{Main Theorem}
\newtheorem{Lem}{Lemma}
\newtheorem{Prop}{Proposition}

\theoremstyle{definition}
\newtheorem{Def}{Definition}
\newtheorem{Note}{Note}
\newtheorem{Ex}{Example}

\theoremstyle{remark}
\newtheorem{notation}{Notation}
\renewcommand{\thenotation}{}

\errorcontextlines=0
\numberwithin{equation}{section}
\renewcommand{\rm}{\normalshape}%

\title[Classical Curvature Flows]%
   {Parabolic Classical Curvature Flows}

\author{Brendan Guilfoyle}
\address{Brendan Guilfoyle\\
          School of Science, Technology, Engineering and Mathematics  \\
          Institute of Technology, Tralee \\
          Clash \\
          Tralee  \\
          Co. Kerry \\
          Ireland.}
\email{brendan.guilfoyle@@ittralee.ie}

\author{Wilhelm Klingenberg}
\address{Wilhelm Klingenberg\\
 Department of Mathematical Sciences\\
 University of Durham\\
 Durham DH1 3LE\\
 United Kingdom.}
\email{wilhelm.klingenberg@@durham.ac.uk }


\date{18 June 2015}

\begin{abstract}
We consider classical curvature flows: 1-parameter families of convex embeddings of the 2-sphere into Euclidean 3-space which evolve by an arbitrary
(non-homogeneous) function of the radii of curvature. The associated flow of the radii of curvature is a second order system of partial differential 
equations which we show decouples to highest order. We determine the conditions for this system to be parabolic and investigate the lower order terms. 
The zeroth order terms are shown to form a Hamiltonian system, which is therefore completely integrable. 

We find conditions on parabolic curvature flows that ensure boundedness of various geometric quantities 
and investigate some examples, including powers of mean curvature flow, Gauss curvature flow and mean radius of curvature flow, as well as 
the non-homogeneous Bloore flow. 

As a new tool we introduce the Radii of Curvature diagram of a surface and its canonical hyperbolic metric. The relationship between these and the
properties of Weingarten surfaces are also discussed. 

\end{abstract}
\thanks{A video of a talk explaining the methods and results of this paper can be found at the following link: \url{https://youtu.be/Sx1T6legtgQ}}
\maketitle
\tableofcontents
\section{Introduction and Results}

The flow of a submanifold embedded in Euclidean space by a specified functional of the eigenvalues of its second fundamental form 
has been studied for decades. For example, there are 
many longtime existence and convergence results for flows by homogeneous symmetric functions of the principal curvatures, such as mean curvature,
inverse mean curvature and Gauss curvature flows \cite{Andrews1} \cite{Chow} \cite{Gerhardt} \cite{gak7} \cite{Huisken} \cite{Schnurer} \cite{Schulze} 
\cite{Smoczyk} \cite{Tso}.

The main concern of this paper is the qualitative study of general non-homogeneous curvature flows for closed convex surfaces in 3-dimensional 
Euclidean space - the domain of classical surface theory. Such equations arise in many physical 
situations, for example, when considering the erosion of a pebble under various types of abrasion \cite{Bloore} \cite{DomGib} \cite{Firey}. 

Define a {\it classical curvature flow} to be a 1-parameter family of smooth convex embeddings $\vec{X}:S^2\times[0,t_1)\rightarrow {\mathbb R}^3$  
such that
\[
\frac{\partial \vec{X}}{\partial t}^\perp=-{\cal K}(\psi,|\sigma|)\;\vec{n},
\]
\[
\vec{X}(S^2,0)=S_0,
\]
where ${\cal K}$ is a given function of $\psi={\textstyle{\frac{1}{2}}}(r_1+r_2)$ and $|\sigma|={\textstyle{\frac{1}{2}}}(r_1-r_2)$, $r_1\geq r_2$
being the radii of curvature of $S=\vec{X}(S^2)$, $\vec{n}$ is the unit normal to $S$ and $S_0$ is an initial convex surface. 

The reasons for this particular combination of radii of curvature will 
become apparent later. A flow is said to be {\it contracting} if ${\cal K}\geq0$ and {\it expanding} if ${\cal K}\leq0$.

A stationary solution of a classical curvature flow is a {\it Weingarten surface}, satisfying the equation ${\cal K}(\psi,|\sigma|)=0$. 
While non-round Weingarten spheres exist (e.g. Hopf spheres \cite{Hopf}), there are many results which state that the only Weingarten spheres
satisfying some particular relationship are round \cite{chern1} \cite{chern2} \cite{HandW1} \cite{Hopf_book} \cite{Voss}.

In general, a Weingarten equation is a second order fully-nonlinear partial differential equation for the support function of the surface. 
Thus ellipticity can be defined, as can the equivalent definition of parabolicity for flows \cite{Lieberman}. 

This paper considers how the radii of curvature, or rather their sum and difference $\psi$ and $|\sigma|$, evolve under a classical 
curvature flow. The pair  $(\psi,|\sigma|)$ take values in the upper half-plane, which we call the {\it radii of curvature space}, and one can visualize 
the dynamics of the flow by considering the image of this map, denoted by $f_t:S^2\rightarrow {\mathbb R}^2_+$ and referred to as the 
{\it radii of curvature (RoC) diagram} of the surface.

The flow on radii of curvature space is a second order system of partial differential equations with many beautiful properties.  Firstly, 
as a consequence of the derived Codazzi-Mainardi equations, 
the system decouples to highest order. Secondly, while the second order differential operator depends upon the 
principal foliation of the surface, ellipticity of the operator is independent of it. 

Furthermore, the  associated {\em classical curvature flow ODE},  
which neglects the spatial derivative terms, is a Hamiltonian system with canonical coordinates
$\psi$ and $|\sigma|$. With the aid of these properties we prove the following. 

Let $\{,\}$ be the Poisson bracket and denote differentiation with respect to the canonical coordinates by 
ordered subscripts.

\vspace{0.1in}

\begin{Thm}\label{t:1stest}
Consider a classical curvature flow with induced flow of radii of curvature: $f_t:S^2\rightarrow {\mathbb R}^2_+$. 

If the flow is contracting ${\cal K}\geq0$ and the function satisfies
\begin{itemize}
\item[(i)] Parabolicity: $-{\cal K}_{10}>|{\cal K}_{01}|$,
\item[(ii)] Convexity: $[{\mbox{Hess}}({\cal K})]\geq0$,
\end{itemize}
then for any function ${\cal H}:{\mathbb R}^2_+\rightarrow{\mathbb R}$ satisfying
\begin{itemize}
\item[(a)] Ellipticity: ${\cal H}_{10}\geq|{\cal H}_{01}|$,
\item[(b)] Convexity: $[{\mbox{Hess}}({\cal H})]\geq0$,
\item[(c)] $\{ {\cal H},{\cal K}\}\geq0$,
\end{itemize}
the following a priori estimate holds for ${\cal H}\circ f_t:S^2\rightarrow {\mathbb R}$:
\[
{\cal H}\circ f_t\leq\max_{S^2}{\cal H}\circ f_0.
\]
If, on the other hand, the flow is expanding ${\cal K}\leq0$ and the function satisfies
\begin{itemize}
\item[(i)] Parabolicity: $-{\cal K}_{10}>|{\cal K}_{01}|$,
\item[(ii)] Concavity: $[{\mbox{Hess}}({\cal K})]\leq0$,
\end{itemize}
then for any function ${\cal H}:{\mathbb R}^2_+\rightarrow{\mathbb R}$ satisfying
\begin{itemize}
\item[(a)] Ellipticity: ${\cal H}_{10}\geq|{\cal H}_{01}|$,
\item[(b)] Concavity: $[{\mbox{Hess}}({\cal H})]\leq0$,
\item[(c)] $\{ {\cal H},{\cal K}\}\leq0$,
\end{itemize}
the following a priori estimate holds:
\[
{\cal H}\circ f_t\geq\min_{S^2}{\cal H}\circ f_0.
\]

\end{Thm}

\vspace{0.1in}

This follows from a more general technical Main Theorem proven in Section \ref{s:tme} where conditions (ii), (a) and (b) are replaced by a single
condition. The proof involves computing the flow of an arbitrary function 
${\cal H}$, studying the sign of lower order terms of the flow and applying the parabolic maximum principle.

There are a variety of applications of this result. For example, clearly for ${\cal H}=-{\cal K}$, conditions (i) and (ii) imply conditions (a)
and (b), while condition (c) is automatically satisfied. In fact, we can drop the convexity condition on ${\cal K}$ and prove

\vspace{0.1in}
\begin{Thm}\label{t:2ndest}
For a parabolic classical curvature flow on $[0,t_1)\times S^2$, the following estimate holds:
\[
|{\cal K}(t)|\geq\min_{S^2}|{\cal K}(0)|.
\]
\end{Thm}
\vspace{0.1in}

By appealing to the more general result, one can also relax the convexity assumption on ${\cal K}$ in Theorem 1 
and prove the following bound on the mean radius of curvature:

\vspace{0.1in}
\begin{Thm}\label{t:3rdest}
For a parabolic classical curvature flow with ${\cal K}+|\sigma|{\cal K}_{01}\geq 0$ and ${\cal K}_{01}+|\sigma|{\cal K}_{02}\geq0$
on $[0,t_1)\times S^2$, 
\[
\psi(t)\leq\max_{S^2}\psi(0).
\]
For a parabolic classical curvature flow with ${\cal K}+|\sigma|{\cal K}_{01}\leq 0$ and ${\cal K}_{01}+|\sigma|{\cal K}_{02}\leq0$
on $[0,t_1)\times S^2$, 
\[
\psi(t)\leq\max_{S^2}\psi(0).
\]
\end{Thm}
\vspace{0.1in}

Similarly, the ellipticity assumption on ${\cal H}$ can be relaxed and still obtain such estimates as:

\vspace{0.1in}

\begin{Thm}\label{t:para}
Consider a parabolic classical curvature flow such that $-{\cal K}_{10}>\epsilon\geq0$. If the flow satisfies 
$-{\cal K}_{10}\geq|\sigma||{\cal K}_{20}|$ on $S^2\times [0,t_1)$, then 
\[
|\sigma(t)|\leq\max_{S^2}|\sigma(0)|e^{-\epsilon t}.
\]
\end{Thm}

\vspace{0.1in}

This indicates the curvature flows that tend to round spheres (for which $|\sigma|=0$), so that the
RoC diagram shrinks to a point on the boundary of the upper half-plane. 

As a further application we prove the non-existence of homothetic solitons other than round spheres, for a wide class of flows. A homothetic 
soliton of a classical flow is a surface that evolves by a dilation under the flow. 

We prove:

\vspace{0.1in}
\begin{Thm}\label{t:soliton}
The only homothetic soliton for a contracting parabolic classical curvature flow with ${\cal K}+|\sigma|{\cal K}_{01}\geq 0$ and 
${\cal K}_{01}+|\sigma|{\cal K}_{02}\geq0$ is the evolving round sphere.

Similarly, the only homothetic soliton for an expanding parabolic classical curvature flow with ${\cal K}+|\sigma|{\cal K}_{01}\leq 0$ and 
${\cal K}_{01}+|\sigma|{\cal K}_{02}\leq0$ is the evolving round sphere.

\end{Thm}
\vspace{0.1in}

Turning to specific classes of flows we establish:

\vspace{0.1in}

\begin{Thm}\label{t:ex}
Consider the following flows: 
\begin{itemize}
\item[(i)] positive powers of mean curvature ${\cal K}=H^n$ for $n\neq0$,
\item[(ii)] positive powers of Gauss curvature  ${\cal K}=K^n$ for $n\neq0$,
\item[(iii)] powers of mean radius of curvature  ${\cal K}=\pm K^n H^{-n}$ for $n\neq0$,
\item[(iv)] linear Weingarten flow ${\cal K}=a+2bH+cK$ for $a,b,c$ positive.
\end{itemize}

These flows are all parabolic.  Linear Weingarten flow is convex, as are powers of mean curvature for $n\geq-1$, powers of Gauss curvature for 
$n\geq{\textstyle{\frac{1}{2}}}$ and powers of mean radius of curvature for $n\geq-1$. Powers of mean radius of curvature are
concave for $n\leq-1$. 

For each of these flows we have the following estimate: 
\[
|{\cal K}(t)|\geq\min_{S^2}|{\cal K}(0)|,
\]
while for the first three flows with $n>0$ 
and the last flow for all positive $a$, $b$ and $c$ we have 
\[
\psi(t)\leq\max_{S^2}\psi(0).
\]
For these values the flows do not admit homothetic solitons, 
other than round spheres.

For negative powers of mean radius of curvature we have 
\[
\psi(t)\geq\min_{S^2}\psi(0),
\]
and there are no homothetic solitons, other than round spheres.

\end{Thm}

\vspace{0.1in}

The rest of this paper is organised as follows: the next section contains the background material on convex surfaces and introduces the radii of
curvature diagram. The connection of this work to classical results on Weingarten surfaces is explored and the
significance of the hyperbolic metric on the RoC diagram is also discussed.

Section 3 derives the evolution equations for the radii of curvature and investigates their properties, while 
Section 4 contains the proofs of Theorems \ref{t:1stest} to \ref{t:soliton}. 
Finally, Section 5 looks at examples of curvature flows and contains the proof of Theorem \ref{t:ex}.

\section{Convex Surfaces and Radii of Curvature Diagram}
\subsection{Classical Surface Theory Redux}

Consider a smooth closed convex surface $S\subset{\mathbb R}^3$ given by a map $X:S^2\rightarrow{\mathbb R}^3$. We now outline our approach
to classical surface theory - further details can be found in  \cite{gak4} and references therein.

Let $\xi$ be the standard complex coordinate on $S^2$ and so, since $S$ is convex, we can use the inverse of the Gauss map $S\rightarrow S^2$ to 
make $\xi$ as a local coordinate on $S$. We refer to these as Gauss coordinates and for all local computations that follow we use them exclusively.

Let $r:S^2\rightarrow {\mathbb R}$ be the support function of $S$ and define the complex derivative
\begin{equation}\label{e:eq1}
F={\textstyle{\frac{1}{2}}}(1+\xi\bar{\xi})^2\bar{\partial} r.
\end{equation}
This is a Lagrangian section of the space of oriented lines in ${\mathbb R}^3$, which can be identified with $T{\mathbb S}^2$
endowed with its canonical neutral K\"ahler structure.

The surface $S$ can be reconstructed from the support function and its derivatives by 
$\vec{X}(\xi,\bar{\xi})=(x^1(\xi,\bar{\xi}),x^2(\xi,\bar{\xi}),x^3(\xi,\bar{\xi}))$ where
\begin{equation}\label{e:coord1}
x^1+ix^2=\frac{2(F-\bar{F}\xi^2)+2\xi(1+\xi\bar{\xi})r}{(1+\xi\bar{\xi})^2}
\qquad
x^3=\frac{- 2(F\bar{\xi}+\bar{F}\xi)+(1-\xi^2\bar{\xi}^2)r}{(1+\xi\bar{\xi})^2}.
\end{equation}
Moving up a derivative, label the complex slopes of $F$ by
\begin{equation}\label{e:eq2}
\psi=r+(1+\xi\bar{\xi})^2\partial\left(\frac{F}{(1+\xi\bar{\xi})^2}\right) \qquad\qquad  \sigma=-\partial\bar{F}.
\end{equation}
By its definition and equation (\ref{e:eq1}), $\psi$ is clearly real. The average and difference of the radii of curvature 
$r_1\geq r_2$ of $S$ can be expressed as:
\begin{equation}\label{e:rocdef}
\psi={\textstyle{\frac{1}{2}}}(r_1+r_2)  \qquad\qquad |\sigma|={\textstyle{\frac{1}{2}}}(r_1-r_2).
\end{equation}
The argument of $\sigma=|\sigma|e^{i\phi}$ gives the principal directions 
of $S$ and, as we will see, is the parameter $\phi$ that appears in the 
differential second order partial differential operator $\triangle_\phi$. 

These quantities satisfy the derived Codazzi-Mainardi equations, which in Gauss coordinates are
\begin{equation}\label{e:comain1}
\partial\psi=-(1+\xi\bar{\xi})^2\bar{\partial}\left(\frac{\sigma}{(1+\xi\bar{\xi})^2}\right).
\end{equation}
This can also be written 
\begin{equation}\label{e:comain1a}
i|\sigma|\partial\phi=\partial|\sigma|+e^{i\phi}\bar{\partial}\psi-\frac{2\bar{\xi}|\sigma|}{1+\xi\bar{\xi}}.
\end{equation}
For future use, note that

\vspace{0.1in}
\begin{Prop}
The following identities hold
\begin{align}\label{e:comain2a}
e^{-i\phi}\partial\partial\psi=&-\partial\bar{\partial}|\sigma|-i|\sigma|\partial\bar{\partial}\phi+ie^{-i\phi}\partial\phi\partial\psi
    -ie^{-i\phi}\partial|\sigma|\bar{\partial}\phi\nonumber\\
  & \qquad+\frac{2\xi\partial|\sigma|}{1+\xi\bar{\xi}}+\frac{2|\sigma|}{(1+\xi\bar{\xi})^2},
\end{align}
\begin{align}\label{e:comain2b}
e^{-i\phi}\partial\partial|\sigma|=&-\partial\bar{\partial}\psi+i|\sigma|e^{-i\phi}\partial\partial\phi+2ie^{-i\phi}\partial\phi\partial|\sigma|
    +|\sigma|e^{-i\phi}(\partial\phi)^2\nonumber\\
  & \qquad+\frac{2\bar{\xi}e^{-i\phi}(\partial|\sigma|-i|\sigma|\partial\phi)}{1+\xi\bar{\xi}}-\frac{2\bar{\xi}^2|\sigma|e^{-i\phi}}{(1+\xi\bar{\xi})^2}.
\end{align}
\end{Prop}
\begin{pf}
The first of these equations is a rearrangement of $\partial$ of $e^{-i\phi}$ times equation (\ref{e:comain1}) and the
second comes from $\bar{\partial}$ of equation (\ref{e:comain1}).

\end{pf}
\vspace{0.1in}
\subsection{The RoC Diagram}\label{s:RoCD}

For a convex surface $S$ in ${\mathbb R}^3$, the radii of curvature give a map $(r_1,r_2):S\rightarrow {\mathbb R}^2_+$. In this context we refer
to ${\mathbb R}^2_+$ as the {\it radii of curvature space}. The image of the map 
$f:S^2\rightarrow {\mathbb R}_+^2$ given by $f(\xi,\bar{\xi})=(\psi(\xi,\bar{\xi}),|\sigma(\xi,\bar{\xi})|)$ we call the {\it radii of curvature
diagram} or {\it RoC diagram} of $S$. 
Many geometric properties of $S$ can be read off its RoC diagram, including convexity and the number of umbilic points.

Figure 1 is the RoC diagram of a generic closed convex surface $S$: it is convex with outward pointing normal and so it lies below the 
diagonal in the first quadrant. The umbilic points on $S$ map to the boundary of the upper half-plane since $|\sigma|=0$ at such points (see 
the second of equations (\ref{e:rocdef})). For a convex surface with inward pointing normals, the RoC diagram would lie below the diagonal 
in the second quadrant.

\vspace{0.1in}
\setlength{\epsfxsize}{5.0in}
\begin{center}
{\mbox{\epsfbox{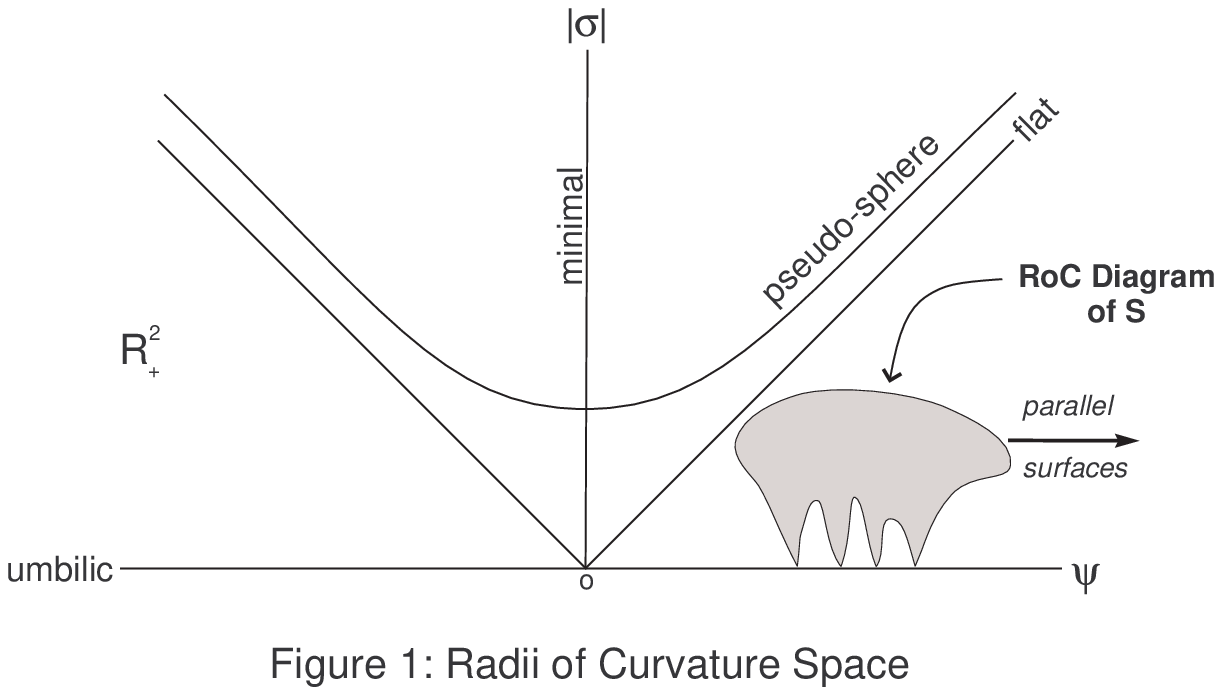}}}
\end{center}
\vspace{0.1in}

\begin{Def}
A surface is {\em Weingarten} if there exists a functional relationship between the radii of curvature: ${\cal K}(\psi,|\sigma|)=0$.
\end{Def}

A convex surface $S$ is Weingarten iff the infinitesimal area of its RoC diagram is zero.
Different types of Weingarten surfaces give rise to different 1-dimensional RoC diagrams, examples of which are also indicated in the figure 
(flat, pseudo-spherical, minimal surfaces). Round spheres have a single point at the boundary as their RoC diagram.

The trio of Weingarten, real analyticity and the behaviour at umbilic points was investigated in the middle of the last century by
a number of authors \cite{chern1} \cite{chern2} \cite{HandW1} \cite{Hopf} \cite{Hopf_book} \cite{Voss} - see also \cite{kuehnel}. 
The curvatures, rather than the radii of curvature, were used in these works and the curvature equivalent of the RoC diagram, called the 
{\it W-diagram}, was introduced. 

The  boundary of the image of $f_t(S)\subset{\mathbb R}^2_+$ is continuous, but in general not smooth at an umbilic point. 
Let $p_0\in S$ be an isolated umbilic point and consider
\[
\kappa(p_0)=\lim_{p\rightarrow p_0}\frac{\psi(p)-\psi(p_0)}{|\sigma(p)|}.
\]
In general, this limit is not well-defined, as it depends upon the direction of approach to the umbilic. As a result, in general the RoC diagram
at an umbilic point is a solid wedge. 

If, however, $S$ is Weingarten, then $\kappa$ is well-defined. Moreover, for real analytic Weingarten surfaces, 
the slope at which the RoC curve strikes the boundary is quantized:

\vspace{0.1in}
\begin{Thm}\cite{Hopf}
Let $S$ be a real analytic Weingarten surface. Then if $p_0$ is an isolated umbilic point, $\kappa(p_0)$ takes one of the following 
values:
\[
-1,\qquad 1, \qquad 0, \qquad \left(\frac{m+1}{m}\right)^{\pm 1},
\]
for $m\in{\mathbb N}$.
\end{Thm}  
\vspace{0.1in}

In addition, 

\vspace{0.1in}
\begin{Def}
A Weingarten surface is {\em special Weingarten} if at its umbilic points  
\[
\frac{\partial{\cal K}}{\partial\lambda_1}.\frac{\partial{\cal K}}{\partial\lambda_2}>0,
\]
where $\lambda_1$ and $\lambda_2$ are the curvatures. 

As will be seen in Note \ref{n:ellip} in Section \ref{s:mainest}, this is the condition that the Weingarten relation is elliptic at the umbilic points.
\end{Def}
\vspace{0.1in}

Weingarten spheres often turn out to be round, for example:

\vspace{0.1in}
\begin{Thm}\cite{HandW1}
Let $S$ be a closed special Weingarten surface of genus zero, which is $C^3$-imbedded in Euclidean space. Then $S$ is a round sphere.
\end{Thm}
\vspace{0.1in}

In contrast to such non-existence results, the Hopf spheres are a 2-parameter family of non-spherical Weingarten surfaces and 
therefore satisfy a non-elliptic relationship. 
For certain integer values of these parameters, the solutions are real analytic, while for non-integers they are smooth non-real analytic surfaces
\cite{kuehnel}.

The hyperbolic metric on ${\mathbb R}^2_+$ plays an interesting role in the RoC space.  Given a surface $S$, constant speed motion along the normal 
lines of the surface (to parallel surfaces) induces a translation in ${\mathbb R}^2_+$ parallel to the boundary, which is a hyperbolic isometry. 

In fact, 

\vspace{0.1in}

\begin{Thm}
Let $f:S^2\rightarrow {\mathbb R}^2_+$ be the RoC diagram of $S$ and  $l:S^2\rightarrow T{\mathbb S}^2$ be the map that takes a point on 
the surface to its oriented normal line, considered as a surface in the space of all oriented lines, $T{\mathbb S}^2$. 

Let $d_{{\mathbb H}^2}A$ be the hyperbolic area on ${\mathbb R}^2_+$ and ${\mathbb G}$ be the canonical invariant neutral metric on $T{\mathbb S}^2$. 
Denote the curvature 2-form of the Lorentz metric induced on $l(S^2)$ by ${\mathbb G}$ by $\Omega_{\mathbb G}$. 

Then 
\[
l_*\Omega_{\mathbb G}=f_*d_{{\mathbb H}^2}A.
\]
\end{Thm}
\begin{pf}
Compare with the curvature expressions in the proof of Main Theorem 3 of \cite{gak5}. In particular, with a slight shift of notation 
($\sigma\leftrightarrow\sigma_0$ and $\psi\leftrightarrow r+\rho_0$), the curvature of the induced metric ${\mathbb G}$ on the Lagrangian
section is
\[
{\mathbb K}=\frac{(1+\xi\bar{\xi})^2}{8|\sigma|^3}(\partial\psi\bar{\partial}|\sigma|-\bar{\partial}\psi\partial|\sigma|),
\]
and the area form is 
\[
d_{\mathbb G}A=\frac{8|\sigma|d\xi\wedge d\bar{\xi}}{(1+\xi\bar{\xi})^2}.
\]
Thus the curvature 2-form is
\[
l_*\Omega_{\mathbb G}=l_*{\mathbb K}\;d_{\mathbb G}A=f_*\frac{d\psi\wedge d|\sigma|}{|\sigma|^2},
\]
which is the pullback of the hyperbolic area 2-form on RoC space, as claimed.
\end{pf}
\vspace{0.1in}

Geometrizing the RoC space as the hyperbolic plane has the added advantage that it places the umbilic points at infinity. Given that many 
results of classical surface theory hold ``away from umbilic points'', such results can be viewed as holding in the hyperbolic plane. Moreover,
it suggests that properties of umbilics may be amenable to exploration using asymptotic methods of the hyperbolic plane. 

It is also worth noting that the divergence between smoothness and real analyticity at umbilic points is further in evidence in the local
version of the Carath\'eodory conjecture \cite{gak6}.

\vspace{0.2in}

\section{Classical Curvature Flows}\label{s:ccf}

\subsection{Evolution Equations}

Consider a classical curvature flow $\vec{X}:S^2\times[0,t_1)\rightarrow {\mathbb R}^3$ such that
\[
\frac{\partial \vec{X}}{\partial t}^\perp =-{\cal{K}}(\psi,|\sigma|)\;\vec{n},
\]
\[
\vec{X}(S^2,0)=S_0,
\]
where ${\cal K}$ is a given function of the radii of curvature, $\vec{n}$ is the unit normal vector to the flowing surface and $S_0$ is an initial 
convex surface. The flow is {\it contracting} if ${\cal K}\geq0$ everywhere and {\it expanding} if ${\cal K}\leq0$ everywhere. 

\vspace{0.1in}
\begin{Prop} \label{p:suppflow}
The support function of $S$ evolves under a classical curvature flow by
\[
\frac{\partial r}{\partial t} =-{\cal{K}}.
\]
\end{Prop}
\begin{pf}
Differentiating equations (\ref{e:coord1}) in time
\[
\frac{\partial}{\partial t}(x^1+ix^2)=\frac{2}{(1+\xi\bar{\xi})^2}\frac{\partial\eta}{\partial t}
-\frac{2\xi^2}{(1+\xi\bar{\xi})^2}\frac{\partial\bar{\eta}}{\partial t}
+\frac{2\xi}{1+\xi\bar{\xi}}\frac{\partial r}{\partial t},
\]
\[
\frac{\partial}{\partial t}x^3=-\frac{2\bar{\xi}}{(1+\xi\bar{\xi})^2}\frac{\partial\eta}{\partial t}
-\frac{2\xi}{(1+\xi\bar{\xi})^2}\frac{\partial\bar{\eta}}{\partial t}
+\frac{1-\xi\bar{\xi}}{1+\xi\bar{\xi}}\frac{\partial r}{\partial t}.
\]
Projecting onto the normal direction $\vec{n}$
\[
\frac{\partial \vec{X}}{\partial t}^\perp=\frac{\partial r}{\partial t}\vec{n}=-{\cal K}\vec{n}.
\]
This yields the stated flow for the support function.
\end{pf}
\vspace{0.1in}

We now compute the evolution of the functions  $\psi$ and $|\sigma|$. For ${\cal K}={\cal K}(\psi,|\sigma|)$ denote derivatives of ${\cal K}$ with 
respect to its arguments by ordered subscripts, so that 
\[
{\cal{K}}_{10}=\frac{\partial{\cal{K}}}{\partial \psi} \qquad {\cal{K}}_{01}=\frac{\partial{\cal{K}}}{\partial |\sigma|}
\qquad
{\cal{K}}_{11}=\frac{\partial^2{\cal{K}}}{\partial \psi\partial|\sigma|}\qquad etc.
\]

\vspace{0.1in}

\begin{Prop}\label{p:rocflow}
The quantities $\psi$ and $|\sigma|$ flow as follows
\[
\frac{\partial \vec{V}}{\partial t}=\triangle_\phi \vec{V}+{\cal Q}(d\vec{V})+{\cal Z}(\vec{V}),
\]
where $\vec{V}=(\psi,|\sigma|)$,
\[
\triangle_\phi={\textstyle{\frac{1}{2}}}(1+\xi\bar{\xi})^2\left[-{\cal K}_{10}\partial\bar{\partial}
 +{\textstyle{\frac{1}{2}}}{\cal K}_{01}(e^{-i\phi}\partial\partial+e^{i\phi}\bar{\partial}\bar{\partial})\right],
\]
\[
{\cal Q}(d\vec{V})=\left[\begin{matrix}
{\cal Q}_1(d\vec{V}) \\
{\cal Q}_2(d\vec{V})
\end{matrix}
\right]
\qquad\qquad
{\cal Z}\left[\begin{matrix}
\psi \\
|\sigma|
\end{matrix}
\right]=\left[\begin{matrix}
-{\cal K}-|\sigma|{\cal K}_{01} \\
|\sigma|{\cal K}_{10} 
\end{matrix}
\right],
\]
and
\begin{align}
{\cal Q}_1&={\textstyle{\frac{1}{2}}}(1+\xi\bar{\xi})^2\left[-{\textstyle{\frac{1}{|\sigma|}}}({\cal K}_{01}
   +|\sigma|{\cal K}_{20})\partial\psi\bar{\partial}\psi
-{\textstyle{\frac{1}{|\sigma|}}}{\cal K}_{01}\left(e^{-i\phi}\partial\psi\partial|\sigma|
   +e^{i\phi}\bar{\partial}\psi\bar{\partial}|\sigma|\right)\right.\nonumber\\
 &\left.\qquad\qquad\qquad-{\cal K}_{11}(\partial\psi\bar{\partial}|\sigma|+\bar{\partial}\psi\partial|\sigma|)
  -{\textstyle{\frac{1}{|\sigma|}}}({\cal{K}}_{01}+|\sigma|{\cal{K}}_{02})\partial|\sigma|\bar{\partial}|\sigma|\right]\nonumber\\
 & \qquad\qquad\qquad+{\textstyle{\frac{1}{2}}}(1+\xi\bar{\xi}){\cal K}_{01}(\bar{\xi}e^{-i\phi}\partial\psi+\xi e^{i\phi}\bar{\partial}\psi),\nonumber\\
&\nonumber\\
{\cal Q}_2&={\textstyle{\frac{1}{2}}}(1+\xi\bar{\xi})^2\left[{\textstyle{\frac{1}{|\sigma|}}}{\cal K}_{10}\partial\psi\bar{\partial}\psi
   +{\textstyle{\frac{1}{2}}}{\cal K}_{20}(e^{-i\phi}(\partial\psi)^2+e^{-i\phi}(\bar{\partial}\psi)^2)\right.\nonumber\\
 &\qquad\qquad\qquad+{\textstyle{\frac{1}{|\sigma|}}}({\cal K}_{10}+|\sigma|{\cal K}_{11})\left(e^{-i\phi}\partial\psi\partial|\sigma|
   +e^{i\phi}\bar{\partial}\psi\bar{\partial}|\sigma|\right)\nonumber\\
 &\left.\qquad\qquad\qquad +\frac{1}{|\sigma|}{\cal K}_{10}\partial|\sigma|\bar{\partial}|\sigma|
     +{\textstyle{\frac{1}{2}}}{\cal K}_{02}(e^{-i\phi}(\partial|\sigma|)^2+e^{i\phi}(\bar{\partial}|\sigma|)^2)\right] \nonumber\\
 & \qquad\qquad\qquad+{\textstyle{\frac{1}{2}}}(1+\xi\bar{\xi}){\cal K}_{01}(\bar{\xi}e^{-i\phi}\partial|\sigma|+\xi e^{i\phi}\bar{\partial}|\sigma|).
\end{align}
\end{Prop}

\begin{pf}

By equation (\ref{e:eq1}) and Proposition \ref{p:suppflow} we have
\[
\frac{\partial F}{\partial t}= {\textstyle{\frac{1}{2}}}(1+\xi\bar{\xi})^2\bar{\partial}\left(\frac{\partial r}{\partial t}\right)
=-{\textstyle{\frac{1}{2}}}(1+\xi\bar{\xi})^2\bar{\partial}{\cal K}=-{\textstyle{\frac{1}{2}}}(1+\xi\bar{\xi})^2({\cal K}_{10}\bar{\partial}\psi+{\cal K}_{01}\bar{\partial}|\sigma|).
\]
Now using the first definition in (\ref{e:eq2})
\begin{align}
\frac{\partial \psi}{\partial t}&=\frac{\partial r}{\partial t}+\partial\left(\frac{\partial F}{\partial t}\right)
-\frac{2\bar{\xi}}{1+\xi\bar{\xi}}\frac{\partial F}{\partial t}\nonumber\\
&=-{\textstyle{\frac{1}{2}}}(1+\xi\bar{\xi})^2\left[{\cal K}_{10}\partial\bar{\partial}\psi
+{\cal K}_{01}\partial\bar{\partial}|\sigma|+{\cal K}_{20}\partial\psi\bar{\partial}\psi+{\cal K}_{11}(\partial|\sigma|\bar{\partial}\psi+\bar{\partial}|\sigma|\partial\psi)\right.\nonumber\\
&\qquad\qquad\qquad\qquad\left.+{\cal K}_{02}\partial|\sigma|\bar{\partial}|\sigma|\right]-{\cal K}.
\end{align}
Substitute the expression for $\partial\bar{\partial}|\sigma|$ by the one from the real part of the derived Codazzi-Mainardi equation 
(\ref{e:comain2a}) and use equation (\ref{e:comain1a}) to remove all derivatives of $\phi$. The result is:
\begin{align}\label{e:psit}
\frac{\partial}{\partial t}\psi=&{\textstyle{\frac{1}{2}}}(1+\xi\bar{\xi})^2\left[-{\cal K}_{10}\partial\bar{\partial}\psi
 +{\textstyle{\frac{1}{2}}}{\cal K}_{01}(e^{-i\phi}\partial\partial\psi+e^{i\phi}\bar{\partial}\bar{\partial}\psi)
 +{\textstyle{\frac{1}{|\sigma|}}}({\cal K}_{01}+|\sigma|{\cal K}_{20})\partial\psi\bar{\partial}\psi\right.\nonumber\\
 &\qquad\qquad\qquad-{\textstyle{\frac{1+\xi\bar{\xi}}{|\sigma|}}}{\cal K}_{01}\left(e^{-i\phi}\partial\psi\partial\left(\frac{|\sigma|}{1+\xi\bar{\xi}}\right)
   +e^{i\phi}\bar{\partial}\psi\bar{\partial}\left(\frac{|\sigma|}{1+\xi\bar{\xi}}\right)\right)\nonumber\\
 &\left.\qquad\qquad\qquad-{\cal K}_{11}(\partial\psi\bar{\partial}|\sigma|+\bar{\partial}\psi\partial|\sigma|)
  -\left({\cal{K}}_{02}+\frac{{\cal{K}}_{01}}{|\sigma|}\right)\partial|\sigma|\bar{\partial}|\sigma|\right]\nonumber\\
 &\qquad\qquad-{\cal K}-|\sigma|{\cal K}_{01}.
\end{align}
This is identical to the expression  for the flow of $\psi$ in the statement of the Proposition, 
given the definitions of $\triangle_\phi$, ${\cal Q}_1$ and ${\cal Z}$.

Turning to the second definition in (\ref{e:eq2}):
\begin{align}
e^{-i\phi}\frac{\partial \sigma}{\partial t}&=-e^{-i\phi}\partial\left(\frac{\partial \bar{F}}{\partial t}\right)\nonumber\\
&=-{\textstyle{\frac{1}{2}}}(1+\xi\bar{\xi})^2e^{-i\phi}\left[{\cal K}_{10}\partial\partial\psi
+{\cal K}_{01}\partial\partial|\sigma|+{\cal K}_{20}(\partial\psi)^2\right.\nonumber\\
&\qquad\qquad\qquad\qquad\qquad\qquad\left.+2{\cal K}_{11}\partial|\sigma|\partial\psi+{\cal K}_{20}(\partial|\sigma|)^2\right]\nonumber\\
&\qquad\qquad+{\textstyle{\frac{1}{2}}}(1+\xi\bar{\xi})e^{-i\phi}\bar{\xi}({\cal K}_{10}\partial\psi+{\cal K}_{01}\partial|\sigma|).
\end{align}

Now substitute the expression for $\partial\bar{\partial}\psi$ by the one from derived Codazzi-Mainardi equation 
(\ref{e:comain2a}) and use equation (\ref{e:comain1a}) to remove all derivatives of $\phi$. The result, taking the real part, is:

\begin{align}\label{e:st}
\frac{\partial}{\partial t}|\sigma|=&{\textstyle{\frac{1}{2}}}(1+\xi\bar{\xi})^2\left[-{\cal K}_{10}\partial\bar{\partial}|\sigma|
 +{\textstyle{\frac{1}{2}}}{\cal K}_{01}(e^{-i\phi}\partial\partial|\sigma|+e^{i\phi}\bar{\partial}\bar{\partial}|\sigma|)
   +{\textstyle{\frac{1}{|\sigma|}}}{\cal K}_{10}\partial\psi\bar{\partial}\psi\right.\nonumber\\
 &\qquad\qquad\qquad+{\textstyle{\frac{1}{2}}}{\cal K}_{20}(e^{-i\phi}(\partial\psi)^2+e^{-i\phi}(\bar{\partial}\psi)^2)\nonumber\\
 &\qquad\qquad\qquad+{\textstyle{\frac{1}{|\sigma|}}}({\cal K}_{10}+|\sigma|{\cal K}_{11})\left(e^{-i\phi}\partial\psi\partial|\sigma|
   +e^{i\phi}\bar{\partial}\psi\bar{\partial}|\sigma|\right)\nonumber\\
 &\left.\qquad\qquad\qquad +{\textstyle{\frac{1}{2}}}{\cal K}_{02}(e^{-i\phi}(\partial|\sigma|)^2
   +e^{i\phi}(\bar{\partial}|\sigma|)^2)+\frac{{\cal K}_{10}}{|\sigma|}\partial|\sigma|\bar{\partial}|\sigma|\right]\nonumber\\
 &\qquad\qquad+{\textstyle{\frac{1}{2}}}(1+\xi\bar{\xi}){\cal K}_{01}(\bar{\xi}e^{-i\phi}\partial|\sigma|+\xi e^{i\phi}\bar{\partial}|\sigma|)
+|\sigma|{\cal K}_{10},
\end{align}
which is the claimed flow for $|\sigma|$.
\end{pf}
\vspace{0.1in}

\subsection{The Classical Curvature O.D.E.}

By Proposition \ref{p:rocflow}, a classical curvature flow gives rise to a second order system of partial differential equations
on radii of curvature space which decouples to top order. 
Necessary and sufficient conditions for parabolicity arise from the second order terms. If we can furthermore put a sign on the quadratic first order 
terms, we can then compare the evolution with the zeroth order terms. 

That is, we are lead to study the behaviour of the following ODE, which we refer to as the {\it classical curvature flow ODE}:

\begin{equation}\label{e:ode1}
\frac{\partial}{\partial t}\psi=-{\cal K}-|\sigma|{\cal K}_{01}
\qquad\qquad\frac{\partial}{\partial t}|\sigma|=|\sigma|{\cal K}_{10}.
\end{equation}

\vspace{0.1in}
\begin{Thm}\label{t:ode}
The classical curvature flow ODE is a Hamiltonian system, with conserved quantity ${\cal I}=|\sigma|{\cal K}$ and canonical coordinates
$\psi$ and $|\sigma|$:
\end{Thm}
\begin{pf}
This follows from noting that
\[
\frac{\partial}{\partial t}\psi=-\frac{\partial}{\partial|\sigma|}{\cal I}
\qquad\qquad
\frac{\partial}{\partial t}|\sigma|=\frac{\partial}{\partial\psi}{\cal I}.
\]
\end{pf}

Thus, if the flow behaves well, it should converge to this flow. In Section 5 we compute the ODE flow for some well-known parabolic 
classical curvature flows.

Conversely, the curvature ODE yields insight into the 
behaviour of the full flow and may suggest quantities that are conserved. 

\vspace{0.1in}
\begin{Cor}
The flow of any function ${\cal H}\circ f_t$ is of the form
\[
\frac{\partial}{\partial t} {\cal H}\circ f_t=\triangle_\phi {\cal H}+\tilde{\cal Q}(d{\cal H})-\{{\cal H},{\cal I}\},\nonumber
\]
where $\{,\}$ are the Poisson brackets associated with the canonical coordinates $(\psi,|\sigma|)$.
\end{Cor}
\vspace{0.1in}

More detailed expression will be given the Main Theorem below.

\vspace{0.2in}

\section{A Priori Estimates}

\subsection{The Main Estimate}\label{s:tme}\label{s:mainest}
In this section we study the behaviour of parabolic classical curvature flows and extract a priori estimates under mild assumptions on the 
functional ${\cal K}$.

\begin{Prop}\label{p:parab}
A classical curvature flow is parabolic iff
\[
-{\cal K}_{10}>|{\cal K}_{01}|.
\]
\end{Prop}
\begin{pf}
For parabolicity, we compute the symbol of the operator $\triangle_\phi$ as follows. Introduce real variables $\xi=x+iy$ so that
\[
\partial\bar{\partial}={\textstyle{\frac{1}{4}}}(\partial_x^2+\partial_y^2)
\qquad\qquad
\partial{\partial}={\textstyle{\frac{1}{4}}}(\partial_x^2-\partial_y^2-2i\partial_x\partial_y).
\]
Then the symbol of 
\[
\triangle_\phi={\textstyle{\frac{1}{2}}}(1+\xi\bar{\xi})^2[-{\cal K}_{10}\partial\bar{\partial}
   +{\textstyle{\frac{1}{2}}}{\cal K}_{10}(e^{-i\phi}\partial{\partial}+e^{i\phi}\bar{\partial}\bar{\partial}),
\]
written in real coordinates is
\[
P={\textstyle{\frac{1}{2}}}(1+\xi\bar{\xi})^2\left[\begin{matrix}
-{\cal K}_{10}+{\cal K}_{01}\cos\phi &  -{\cal K}_{01}\sin\phi\\
 -{\cal K}_{01}\sin\phi &-{\cal K}_{10}-{\cal K}_{01}\cos\phi 
\end{matrix}
\right].
\]
This is elliptic if $P(X,X)=0$ implies that $X=0$. In other words
\[
\det P={\textstyle{\frac{1}{4}}}(1+\xi\bar{\xi})^4({\cal K}_{10}^2-{\cal K}_{01}^2)\neq0.
\]
In order for the operator ${\textstyle{\frac{\partial}{\partial t}}}-\triangle_\phi$ to be parabolic, we must also require that $-{\cal K}_{10}>0$.
\end{pf}
\vspace{0.1in}

\begin{Note}\label{n:ellip}
This definition of parabolicity agrees with the usual definition of ellipticity for fully non-linear second order partial differential equations.
That is, an equation involving the second derivatives of $r$
\[
{\cal F}(\partial_x^2r,\partial_x\partial_yr,\partial^2_yr)=0,
\]
is {\it elliptic} if ${\cal F}_1{\cal F}_3-{\cal F}_2^2>0$, where a subscript represents differentiation with respect to the arguments of ${\cal F}$.

In our case, we have (see the definitions (\ref{e:eq1}) and (\ref{e:eq2}) of $\psi$ and $|\sigma|$ in terms of $r$ and switching to real variables 
as above)
\[
{\cal F}(\partial_x^2r,\partial_x\partial_yr,\partial^2_yr)={\cal K}\left(\partial_x^2r+\partial^2_yr,\sqrt{\partial_x^2r-\partial^2_yr+(\partial_x\partial_yr)^2}\right)=0,
\]
for which ${\cal F}_1{\cal F}_3-{\cal F}_2^2>0$ is easily found to be equivalent to ${\cal K}_{10}^2-{\cal K}_{01}^2>0$.

As we saw in Section \ref{s:RoCD}, a Weingarten surface is {\it special Weingarten} if at its umbilic points  
\[
\frac{\partial{\cal K}}{\partial\lambda_1}.\frac{\partial{\cal K}}{\partial\lambda_2}>0,
\]
where $\lambda_1$ and $\lambda_2$ are the curvatures. It is not hard to see that, in terms of the canonical coordinates $(\psi,|\sigma|)$,
\[
\frac{\partial{\cal K}}{\partial\lambda_1}.\frac{\partial{\cal K}}{\partial\lambda_2}={\textstyle{\frac{1}{4}}}(\psi^2-|\sigma|^2)
    ({\cal K}_{10}^2-{\cal K}_{01}^2),
\]
and so {\it special Weingarten} for a convex Weingarten surface is equivalent to {\it elliptic at the umbilic points}.
\end{Note}
\vspace{0.1in}

We now prove the main estimate:

\vspace{0.1in}
\noindent{\bf Main Theorem}: 

{\it 
Consider a classical parabolic curvature flow with induced flow of radii of curvature $f_t:S^2\rightarrow {\mathbb R}^2_+$. 
For ${\cal H}:{\mathbb R}^2_+\rightarrow{\mathbb R}$  define
\begin{equation}\label{e:ab1}
A=({\cal H}_{10}^2+{\cal H}_{01}^2)\{{\cal H},{\cal K}\}-|\sigma|({\cal K}_{10}{\mbox{Hess}}_{\cal H}(d{\cal H})
   -{\cal H}_{10}{\mbox{Hess}}_{\cal K}(d{\cal H})),
\end{equation}
\begin{equation}\label{e:ab2}
B=-2{\cal H}_{10}{\cal H}_{01}\{{\cal H},{\cal K}\}+|\sigma|({\cal K}_{01}{\mbox{Hess}}_{\cal H}(d{\cal H})
   -{\cal H}_{01}{\mbox{Hess}}_{\cal K}(d{\cal H})),
\end{equation}
and 
\begin{equation}\label{e:ab3}
{\mbox{Hess}}_{\cal K}(d{\cal H})={\cal K}_{20}{\cal H}_{01}^2-2{\cal K}_{11}{\cal H}_{01}{\cal H}_{10}+{\cal K}_{02}{\cal H}_{10}^2.
\end{equation}
Suppose that ${\cal H}$ satisfies 
\begin{itemize}
\item[(1a)]\hspace{0.2in} $\{ {\cal H},{\cal I}\}\geq0$,
\item[(2a)]\hspace{0.2in} $A\geq|B|$.
\end{itemize}
Then the following a priori estimate holds:
\[
{\cal H}\circ f_t\leq\max_{S^2}{\cal H}\circ f_0.
\]
On the other hand if ${\cal H}$ satisfies 
\begin{itemize}
\item[(1b)]\hspace{0.2in} $\{ {\cal H},{\cal I}\}\leq0$,
\item[(2b)]\hspace{0.2in} $A\leq-|B|$.
\end{itemize}
then 
\[
{\cal H}\circ f_t\geq\min_{S^2}{\cal H}\circ f_0.
\]
}
\begin{pf}

The flow of the function ${\cal H}\circ f_t$ is computed as follows.
\begin{align}
\frac{\partial}{\partial t} {\cal H}\circ f_t&={\cal H}_{10}\frac{\partial}{\partial t} {\psi}+{\cal H}_{01}\frac{\partial}{\partial t} {|\sigma|}\nonumber\\
&={\cal H}_{10}(\triangle_\phi \psi+{\cal Q}_1-{\cal K}-|\sigma|{\cal K}_{01}))+{\cal H}_{01}(\triangle_\phi |\sigma|+{\cal Q}_2+|\sigma|{\cal K}_{10})\nonumber\\
&={\cal H}_{10}(\triangle_\phi \psi+{\cal Q}_1)+{\cal H}_{01}(\triangle_\phi |\sigma|+{\cal Q}_2)-\{{\cal H},{\cal I}\}\nonumber\\
&=\triangle_\phi {\cal H}+{\cal Q}_3(d{\cal H})-{\cal Z}_2-\{{\cal H},{\cal I}\},\nonumber
\end{align}
where on the second line we have used the flow equations in Proposition 5. Here ${\cal Q}_3(d{\cal H})=0$ if $d{\cal H}=0$ and  
\[
{\cal Z}_2=\frac{(1+\xi\bar{\xi})^2}{2|\sigma|{\cal H}_{10}^2}\left[A\partial|\sigma|\bar{\partial}|\sigma|+{\textstyle{\frac{1}{2}}}B(e^{-i\phi}(\partial|\sigma|)^2+e^{i\phi}(\bar{\partial}|\sigma|)^2)\right],
\]
with $A$ and $B$ given by expressions (\ref{e:ab1}) to (\ref{e:ab3}).

Thus, at the maximum or minimum value of ${\cal H}$,
\[
\left(\frac{\partial}{\partial t} -\triangle_\phi\right){\cal H}\circ f_t=-{\cal Z}_2-\{{\cal H},{\cal I}\}.
\]

If $A\geq|B|$ then ${\cal Z}_2\geq0$ and if in addition $\{{\cal H},{\cal I}\}\geq0$, the estimate follows by the parabolic maximum principle.

If $A\leq-|B|$ then ${\cal Z}_2\leq0$ and if $\{{\cal H},{\cal I}\}\leq0$, the estimate follows.

\end{pf}
\vspace{0.1in}

\subsection{Applications}

More geometric assumptions can be used to obtain a priori bounds:

\vspace{0.1in}
\noindent{\bf Theorem 1}: 

{\it 
Consider a classical curvature flow with induced flow of radii of curvature: $f_t:S^2\rightarrow {\mathbb R}^2_+$. 

If the flow is contracting (${\cal K}\geq0$) and the function satisfies
\begin{itemize}
\item[(a)] Parabolicity: $-{\cal K}_{10}>|{\cal K}_{01}|$,
\item[(b)] Convexity: $[{\mbox{Hess}}({\cal K})]\geq0$,
\end{itemize}
then for any function ${\cal H}:{\mathbb R}^2_+\rightarrow{\mathbb R}$ satisfying
\begin{itemize}
\item[(i)] Ellipticity: ${\cal H}_{10}\geq|{\cal H}_{01}|$,
\item[(ii)] Convexity: $[{\mbox{Hess}}({\cal H})]\geq0$,
\item[(iii)] Poisson: $\{ {\cal H},{\cal K}\}\geq0$,
\end{itemize}
the following a priori estimate holds for ${\cal H}\circ f_t:S^2\rightarrow {\mathbb R}$:
\[
{\cal H}\circ f_t\leq\max_{S^2}{\cal H}\circ f_0.
\]
If, on the other hand, the flow is expanding (${\cal K}\leq0$) and the function satisfies
\begin{itemize}
\item[(a)] Parabolicity: $-{\cal K}_{10}>|{\cal K}_{01}|$,
\item[(b)] Concavity: $[{\mbox{Hess}}({\cal K})]\leq0$,
\end{itemize}
then for any function ${\cal H}:{\mathbb R}^2_+\rightarrow{\mathbb R}$ satisfying
\begin{itemize}
\item[(i)] Ellipticity: ${\cal H}_{10}\geq|{\cal H}_{01}|$,
\item[(ii)] Concavity: $[{\mbox{Hess}}({\cal H})]\leq0$,
\item[(iii)] Poisson: $\{ {\cal H},{\cal K}\}\leq0$,
\end{itemize}
the following a priori estimate holds:
\[
{\cal H}\circ f_t\geq\min_{S^2}{\cal H}\circ f_0.
\]
}
\begin{pf}
Assume that the flow is contracting, so that ${\cal K}\geq0$ and assume that the conditions (a), (b), (i), (ii) and (iii) of the Theorem's statement hold.
Then, by conditions (i) and (iii) 
\[
\{{\cal H},{\cal I}\}={\cal H}_{10}{\cal K}+|\sigma|\{{\cal H},{\cal K}\}\geq0,
\]
so that condition (1a) in the Main Theorem holds.

To see that condition (2a) also holds, compare the expressions for $A$ and $B$ in equations (\ref{e:ab1}) and (\ref{e:ab2}) term by term. Note that: the 
first term of $A$ is positive by condition (iii) and 
dominates the first term of $B$. By condition (a), the second term of $A$ is also positive and dominates the second term of $B$. 

By conditions (a) and (ii), the third term of $A$ is positive and dominates the third term of $B$. The final term in $A$ is positive and 
dominates the final term of $B$ due to conditions (i) and (b). 

Thus conditions (1a) and (2a) of the Main Theorem hold, and we can apply it to yield the stated result.

The proof of the expanding case is analogous with opposite inequalities.
\end{pf}
\vspace{0.1in}

As a consequence we have

\vspace{0.1in}
\noindent{\bf Theorem 2}:

{\it
For a parabolic classical curvature flow on $[0,t_1)\times S^2$, the following estimate holds:
\[
|{\cal K}(t)|\geq\min_{S^2}|{\cal K}(0)|.
\]
}
\vspace{0.1in}
\begin{pf}
Let ${\cal H}=-{\cal K}$, then compute that $A=B=0$ and $\{{\cal H},{\cal I}\}=-{\cal K}{\cal K}_{10}$. Thus for a contracting parabolic flow  
$\{{\cal H},{\cal I}\}\geq0$, so conditions (1a) and (2a) hold  and applying the Main Theorem we obtain ${\cal K}(t)\geq\min_{S^2}{\cal K}(0)$. 

In the expanding parabolic case $\{{\cal H},{\cal I}\}\leq0$ and so (1b) and (2b) hold and again apply the Main Theorem.

\end{pf}
\vspace{0.1in}

Note that this lower bound is a tautology unless ${\cal K}$ has a fixed sign. 

We also have:

\vspace{0.1in}
\noindent{\bf Theorem 3}:

{\it
For a parabolic classical curvature flow with ${\cal K}+|\sigma|{\cal K}_{01}\geq 0$ and ${\cal K}_{01}+|\sigma|{\cal K}_{02}\geq0$
on $[0,t_1)\times S^2$
\[
\psi(t)\leq\max_{S^2}\psi(0).
\]

For a parabolic classical curvature flow with ${\cal K}+|\sigma|{\cal K}_{01}\leq 0$ and ${\cal K}_{01}+|\sigma|{\cal K}_{02}\leq0$
on $[0,t_1)\times S^2$
\[
\psi(t)\geq\min_{S^2}\psi(0).
\]

}
\begin{pf}
This follows from the Main Theorem by setting ${\cal H}=\psi$ and noting that $A={\cal K}_{01}+|\sigma|{\cal K}_{02}$, $B=0$ and 
$\{{\cal H},{\cal I}\}={\cal K}+|\sigma|{\cal K}_{01}$ so that conditions (1a) and (2a) hold for contracting flows under the stated conditions, 
while (1b) and (2b) hold for expanding flows.
\end{pf}
\vspace{0.1in}

An upper bound on the deviation from roundness can also be found:

\vspace{0.1in}
\noindent{\bf Theorem 4}:

{\it
For a parabolic curvature flow with $-{\cal K}_{10}>\epsilon\geq0$ and $-{\cal K}_{10}\geq|\sigma||{\cal K}_{20}|$ on $[0,t_1)\times S^2$, we have 
\[
|\sigma(t)|\leq\max_{S^2}|\sigma(0)|e^{-\epsilon t}.
\]
} 
\begin{pf}
This follows from the Main Theorem by setting ${\cal H}=|\sigma|$ and noting that $A=-{\cal K}_{10}$, $B=-|\sigma|{\cal K}_{20}$ and 
$\{{\cal H},{\cal I}\}=-|\sigma|{\cal K}_{10}$.
\end{pf}
\vspace{0.1in}

Thus, strictly parabolic flows that satisfy ${\cal K}_{10}^2\geq|\sigma|^2{\cal K}_{20}^2$ tend to umbilicity. For expanding flows this would be a
plane at infinity, while for contracting flows, this would be a round sphere.

\subsection{Non-existence of Homothetic Solitons}

A homothetic soliton is a surface $S$ that satisfies the equation
\begin{equation}\label{e:soliton}
\lambda r={\cal K},
\end{equation}
for some constant $\lambda$, which is  positive if the flow is contracting and negative if the flow is expanding, such that if we flow $S$ by 
${\cal K}$, then it simply scales the surface about the origin. 

\vspace{0.1in}
\begin{Thm}\label{t:round}
If at $t=0$ the surface $S$ is a round sphere, then for as long as the classical curvature flow exists, it remains a round sphere with radius $R(t)$ 
evolving by
\[
\frac{dR}{dt}=-{\cal K}(R,0). 
\]
\end{Thm}
\begin{pf}
It is well-known that the only convex totally umbilic surface is a round sphere. That is, $S$ is umbilic iff $|\sigma|=0$, which by equation
(\ref{e:comain1a}) implies that $\psi=constant=R$.

By a translation, let the centre of the initial sphere lie at the origin, so that $F(0)=0$ and $r(0)=R_0$. By the evolution equations
we see that
\[
\frac{\partial F}{\partial t}(0)=0,
\]
and so, by uniqueness of solutions to ODE's, it remains a round sphere centred at the origin, $r=R(t)$ and the flow of the radius is as stated.
\end{pf}
\vspace{0.1in}

We now prove a non-existence results for homothetic solitons.

\vspace{0.1in}
\noindent{\bf Theorem 5}:

{\it 
The only homothetic soliton for a contracting parabolic classical curvature flow with ${\cal K}+|\sigma|{\cal K}_{01}\geq 0$ and 
${\cal K}_{01}+|\sigma|{\cal K}_{02}\geq0$ is the evolving round sphere given in Theorem \ref{t:round}.

Similarly, the only homothetic soliton for an expanding parabolic classical curvature flow with ${\cal K}+|\sigma|{\cal K}_{01}\leq 0$ and 
${\cal K}_{01}+|\sigma|{\cal K}_{02}\leq0$ is the evolving round sphere.
}
\vspace{0.1in}
\begin{pf}
Differentiating equation (\ref{e:soliton}) in a manner similar to the computations in Proposition \ref{p:rocflow} we find that
\begin{align}\label{e:psit}
\lambda\psi=&{\textstyle{\frac{1}{2}}}(1+\xi\bar{\xi})^2\left[-{\cal K}_{10}\partial\bar{\partial}\psi
 +{\textstyle{\frac{1}{2}}}{\cal K}_{01}(e^{-i\phi}\partial\partial\psi+e^{i\phi}\bar{\partial}\bar{\partial}\psi)
 +{\textstyle{\frac{1}{|\sigma|}}}({\cal K}_{01}+|\sigma|{\cal K}_{20})\partial\psi\bar{\partial}\psi\right.\nonumber\\
 &\qquad\qquad\qquad-{\textstyle{\frac{1+\xi\bar{\xi}}{|\sigma|}}}{\cal K}_{01}\left(e^{-i\phi}\partial\psi\partial\left(\frac{|\sigma|}{1+\xi\bar{\xi}}\right)
   +e^{i\phi}\bar{\partial}\psi\bar{\partial}\left(\frac{|\sigma|}{1+\xi\bar{\xi}}\right)\right)\nonumber\\
 &\left.\qquad\qquad\qquad-{\cal K}_{11}(\partial\psi\bar{\partial}|\sigma|+\bar{\partial}\psi\partial|\sigma|)
  -\left({\cal{K}}_{02}+\frac{{\cal{K}}_{01}}{|\sigma|}\right)\partial|\sigma|\bar{\partial}|\sigma|\right]\nonumber\\
 &\qquad\qquad-{\cal K}-|\sigma|{\cal K}_{01}\nonumber.
\end{align}
At the maximum and minimum value of $\psi$ we therefore have
\[
\lambda\psi=\triangle_\phi\psi-{\textstyle{\frac{1}{2}}}(1+\xi\bar{\xi})^2
   \left({\cal{K}}_{02}+\frac{{\cal{K}}_{01}}{|\sigma|}\right)\partial|\sigma|\bar{\partial}|\sigma|-{\cal K}-|\sigma|{\cal K}_{01},
\]
so, under the assumptions of the Theorem, at a maximum of a contracting flow or a minimum of an expanding flow $\psi\leq0$, which is impossible.
\end{pf}
\vspace{0.2in}

\section{Examples}

In this section we consider the following classical curvature flows: powers of mean curvature flow 
\[
{\cal K}=\pm H^n=\pm\left(\frac{r_1+r_2}{r_1r_2}\right)^n=\pm\frac{\psi^n}{(\psi^2-|\sigma|^2)^n},
\]
powers of Gauss curvature flow
\[
{\cal K}=\pm K^n=\pm\frac{1}{(r_1r_2)^n}=\pm\frac{1}{(\psi^2-|\sigma|^2)^n},
\]
power of mean radius of curvature flow
\[
{\cal K}=\pm\left(\frac{H}{K}\right)^n=\pm(r_1+r_2)^n=\pm\psi^n,
\]
and linear Weingarten flow
\[
{\cal K}=a+2bH+cK=a+\frac{c+2b(r_1+r_2)}{r_1r_2}=a+\frac{2b\psi+c}{\psi^2-|\sigma|^2}
\]
where we take the positive (negative) sign on the first three flows for $n>0$ ($n<0$) respectively, and $a,b,c$ are positive constants. 
In fact, the linear Weingarten flow is used as an abrasion model under the assumptions that $a=1$ and $b^2>c>0$,
when it is called the Bloore flow \cite{Bloore} \cite{DomGib} \cite{Galvez}.

\vspace{0.1in}
\noindent{\bf Theorem 6}:

{\it
Consider the above flows: powers of mean curvature, Gauss curvature, powers of mean radius of curvature and the linear Weingarten flow.

These flows are all parabolic.  Linear Weingarten flow is convex, as are powers of mean curvature for $n\geq-1$, powers of Gauss curvature for 
$n\geq{\textstyle{\frac{1}{2}}}$ and powers of mean radius of curvature for $n\geq-1$. Powers of mean radius of curvature are
concave for $n\leq-1$. 

For each of these flows we have the following estimate: 
\[
|{\cal K}(t)|\geq\min_{S^2}|{\cal K}(0)|,
\]
while for the first three flows with $n>0$ 
and the last flow for all positive $a$, $b$ and $c$ we have 
\[
\psi(t)\leq\max_{S^2}\psi(0).
\]
For these values the flows do not admit homothetic solitons, 
other than round spheres.

For negative powers of mean radius of curvature we have 
\[
\psi(t)\geq\min_{S^2}\psi(0),
\]
and there are no homothetic solitons, other than round spheres.
}
\begin{pf}
By direct computation we find the first derivatives as reported in the first Table of the Appendix. In each case note that
$-{\cal K}_{10}>|{\cal K}_{01}|$ and so, as long as the surface remains convex, the flows are parabolic. Thus the first estimate follows 
from Theorem 2.

Moving to second derivatives, the results for each flow are given in the second and third Tables of the Appendix.
Convexity for the stated values of $n$ follows from checking that $|{\mbox{Hess }}{\cal K}|\geq0$ and ${\cal K}_{20}\geq0$ for each flow.

Finally the second estimate and the non-existence of homothetic solutions follows from Theorems 3 and 5 respectively as long as the conditions
${\cal K}+|\sigma|{\cal K}_{01}\geq0$ and ${\cal K}_{01}+|\sigma|{\cal K}_{02}\geq0$ for the contracting flows and 
${\cal K}+|\sigma|{\cal K}_{01}\leq0$ and ${\cal K}_{01}+|\sigma|{\cal K}_{02}\leq0$ for the expanding flows.

\end{pf}
\vspace{0.1in}

\subsection{The ODE Flow}

By Theorem \ref{t:ode}, the classical curvature ODE is a Hamiltonian system with conserved quantity ${\cal I}=|\sigma|{\cal K}$. Thus the
flowlines of the ODE are given by ${\cal I}=constant$.

Figure 2 contains plots of these flowlines in RoC space for mean curvature flow ($n=1$), Gauss curvature flow ($n=1$), mean radius of curvature
flow ($n=-1$) and linear Weingarten flow ($a=1$, $b=2$ and $c=1$).  

The flowlines indicate contraction to a point for the first two and the last flows, and expansion to a sphere at infinity for the mean radius of
curvature flow. While the full PDE will deviate from these flowlines by amounts determined by the spatial derivatives of the radii of curvature,
the qualitative behaviour of the flows are evident in the ODE.

\vspace{0.1in}
\setlength{\epsfxsize}{5.0in}
\begin{center}
{\mbox{\epsfbox{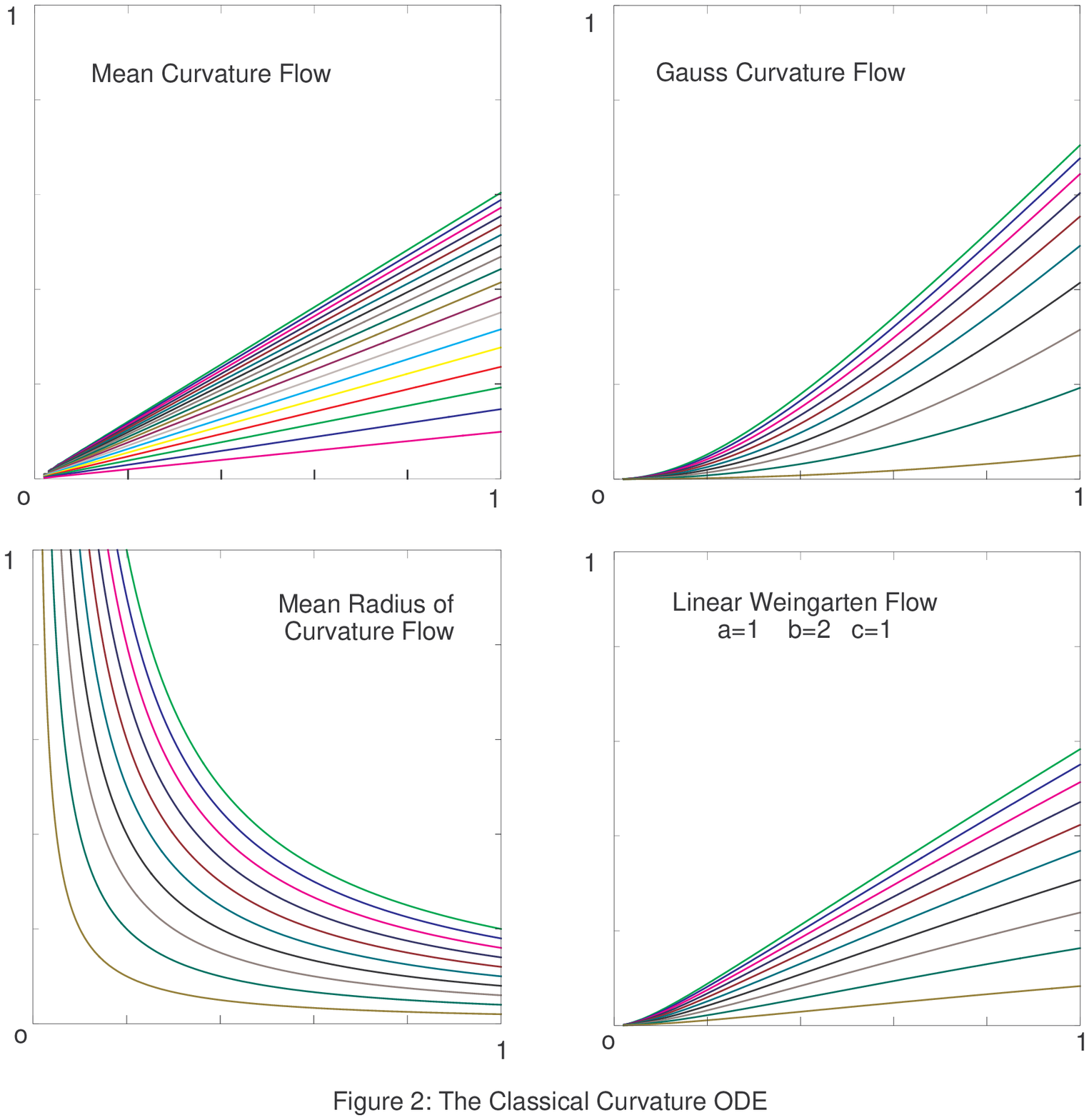}}}
\end{center}
\vspace{0.1in}

\section*{Appendix}

\vspace{0.2in}
\begin{tabular}{|c||c|c|c|c|}
\hline
&&&&\\
& Power of & Power of &Power of Mean & Linear \\
&Mean Curvature & Gauss Curvature &Radii of Curvature&Weingarten\\
&&&&\\
\hline
\hline
&&&&\\
${\cal K}$& $\pm\frac{\psi^n}{(\psi^2-|\sigma|^2)^n}$&$ \pm\frac{1}{(\psi^2-|\sigma|^2)^n}$&$\pm\frac{1}{\psi^n} $&$ a+\frac{2b\psi+c}{\psi^2-|\sigma|^2}$\\
&&&&\\
\hline
&&&&\\
${\cal K}_{10}$& $-\frac{|n|\psi^{n-1}(\psi^2+|\sigma|^2)}{(\psi^2-|\sigma|^2)^{n+1}}$ & $-\frac{2|n|\psi}{(\psi^2-|\sigma|^2)^{n+1}}$ & $-\frac{|n|}{\psi^{n+1}}$ & $-\frac{2(b(\psi^2+|\sigma|^2)+c\psi)}{(\psi^2-|\sigma|^2)^2}$ \\
&&&&\\
\hline
&&&&\\
${\cal K}_{01}$& $\frac{2|n||\sigma|\psi^n}{(\psi^2-|\sigma|^2)^{n+1}}$ & $\frac{2|n||\sigma|}{(\psi^2-|\sigma|^2)^{n+1}}$ & 0 &  $\frac{2|\sigma|(c+2b\psi)}{(\psi^2-|\sigma|^2)^2}$\\
&&&&\\
\hline
\end{tabular}
\vspace{0.2in}

\begin{tabular}{|c||c|c|}
\hline
&&\\
& Power of & Power of  \\
&Mean Curvature & Gauss Curvature \\
&&\\
\hline
\hline
&&\\
${\cal K}_{11}$& $-\frac{2|n||\sigma|\psi^{n-1}[(n+2)\psi^2+n|\sigma|^2]}{(\psi^2-|\sigma|^2)^{n+2}}$
&$ -\frac{4|n|(n+1)|\sigma|}{(\psi^2-|\sigma|^2)^{n+2}}$\\
&&\\
\hline
&&\\
${\cal K}_{20}$& $\frac{|n|\psi^{n-2}[(n+1)\psi^4+2(n+2)\psi^2|\sigma|^2+(n-1)|\sigma|^4]}{(\psi^2-|\sigma|^2)^{n+2}}$ & $\frac{2|n|[(2n+1)\psi^2+|\sigma|^2]}{(\psi^2-|\sigma|^2)^{n+2}}$  \\
&&\\
\hline
&&\\
${\cal K}_{02}$& $\frac{2|n|\psi^n(\psi^2+(2n+1)|\sigma|^2)}{(\psi^2-|\sigma|^2)^{n+2}}$ & $\frac{2|n|[\psi^2+(2n+1)|\sigma|^2]}{(\psi^2-|\sigma|^2)^{n+2}}$ \\
&&\\
\hline
&&\\
$|{\mbox{Hess }}{\cal K}|$& $\frac{2n^2(n+1)\psi^{2n-2}}{(\psi^2-|\sigma|^2)^{2n+1}}$ & $\frac{4n^2(2n+1)}{(\psi^2-|\sigma|^2)^{2n+2}}$\\
&&\\
\hline
\end{tabular}
\vspace{0.2in}

\begin{tabular}{|c||c|c|}
\hline
&&\\
 &Power of Mean & Linear \\
 &Radii of Curvature&Weingarten\\
&&\\
\hline
\hline
&&\\
${\cal K}_{11}$&0 &$ -\frac{4|\sigma|[b(3\psi^2+|\sigma|^2)+2c\psi]}{(\psi^2-|\sigma|^2)^3}$\\
&&\\
\hline
&&\\
${\cal K}_{20}$ & $\frac{|n|(n+1)}{\psi^{n+2}}$ & $\frac{2[2b\psi(\psi^2+3|\sigma|^2)+c(3\psi^2+|\sigma|^2)]}{(\psi^2-|\sigma|^2)^3}$ \\
&&\\
\hline
&&\\
${\cal K}_{02}$& 0 &  $\frac{2[2b\psi(\psi^2+3|\sigma|^2)+c(\psi^2+3|\sigma|^2)]}{(\psi^2-|\sigma|^2)^3}$\\
&&\\
\hline
&&\\
$|{\mbox{Hess }}{\cal K}|$& 0 &  $\frac{4[4b^2(\psi^2-|\sigma|^2)+8bc\psi|\sigma|+3c^2]}{(\psi^2-|\sigma|^2)^4}$\\
&&\\
\hline
\end{tabular}

\end{document}